\documentclass{amsart}
\usepackage{amssymb}
\usepackage{a4}
\begin{document}
%\date{\version}
\newtheorem{theorem}{Theorem}[section]
\newtheorem{lemma}[theorem]{Lemma}
\newtheorem{remark}[theorem]{Remark}
\newtheorem{definition}[theorem]{Definition}
\newtheorem{corollary}[theorem]{Corollary}
\newtheorem{example}[theorem]{Example}
\font\pbglie=eufm10
\def\SS{\mathcal{S}}
\def\RRR{\text{\pbglie R}}
\def\Rank{\operatorname{Rank}}
\def\Pspan{\operatorname{Span}}
\def\Range{\operatorname{Range}}
\makeatletter
  \renewcommand{\theequation}{%
   \thesection.\alph{equation}}
  \@addtoreset{equation}{section}
 \makeatother
\title[Spacelike Jordan Osserman Manifolds]
{Nilpotent Spacelike Jorden Osserman pseudo-Riemannian manifolds}
\author{P. Gilkey and S. Nik\v cevi\'c}
\begin{address}{PG: Mathematics Department, University of Oregon,
Eugene Or 97403 USA.\newline Email: {\it gilkey@darkwing.uoregon.edu}}
\end{address}
\begin{address}{SN: Mathematical Institute, Sanu,
Knez Mihailova 35, p.p. 367,
11001 Belgrade,
Yugoslavia 
%stana@sfb288.math.tu-berlin.de
%enikcevi@ubbg.etf.bg.ac.yu
\newline Email: {\it stanan@mi.sanu.ac.yu}}\end{address}
\begin{abstract}Pseudo-Riemannian manifolds of balanced signature which are both spacelike
and timelike Jordan  Osserman nilpotent of order $2$ and of order $3$ have been
constructed previously. In this short note, we shall construct pseudo-Riemannian manifolds
of signature $(2s,s)$ for any $s\ge2$ which are spacelike Jordan Osserman
nilpotent of order 3 but which are not timelike Jordan Osserman. Our example and techniques
are quite different from known previously both in  that they are not in neutral signature and
that the manifolds constructed will be spacelike but not timelike Jordan 
Osserman.\end{abstract}
\maketitle

\section{Introduction} 

Let $(M,g)$ be a pseudo-Riemannian manifold of signature $(p,q)$. Let
$$S^\pm(M,g):=\{x\in TM:(x,x)=\pm1\}$$
be the bundles of unit spacelike and unit timelike vectors, respectively.  Let $R$ be the
associated Riemann curvature tensor. If $x\in T_PM$, then the {\it Jacobi operator}
$J(x)$ is the self-adjoint linear map of $T_PM$ which is characterized by the identity:
\begin{equation}g(J(x)y,z)=R(y,x,x,z).\label{eqn-1.a}\end{equation}
One says that $(M,g)$ is {\it spacelike Osserman} or {\it
timelike Osserman} if the eigenvalues of $J$ are constant on $S^+(M,g)$ or on $S^-(M,g)$,
respectively. These are equivalent notions if $p\ge1$ and $q\ge1$ \cite{refGil} so such
manifolds are simply said to be {\it Osserman}. 

If $p=0$, and similarly if $q=0$, then one is in the {\it Riemannian
setting}. If $(M,g)$ is a rank $1$ symmetric space or if $(M,g)$ is flat, then the local
isometries of
$(M,g)$ act transitively on $S^+(M,g)$ so the eigenvalues of $J$ are constant on
$S^+(M,g)$. Osserman \cite{refOss} wondered if the converse
held. Work of Chi \cite{refChi} and of
Nikolayevsky \cite{refNik} has shown this to be the case if the dimension  is different from $8$ and $16$.

If $p=1$, and similarly if $q=1$, then one is in the
{\it Lorentzian setting}. Bla\v zi\'c, Bokan and Gilkey \cite{refBBG} and 
Garc\'{\i}a--R\'{\i}o, Kupeli and V\'azquez-Abal \cite{refGKV} have shown that Lorentzian
Osserman manifolds have constant sectional curvature.

The situation is quite different in
the higher signature setting where $p\ge2$ and $q\ge2$. There exist Osserman
pseudo-Riemannian manifolds which are not symmetric spaces
\cite{refBBGZ,refBV03,refBCG,refBCGHV,refGVV}; we refer to \cite{refGRKVL} for an excellent
and quite comprehensive treatment of the subject.

In the higher signature setting, it is natural to impose a more restrictive hypothesis and
study the Jordan normal form of the Jacobi operator. We say that $(M,g)$ is {\it spacelike
Jordan Osserman} or is {\it timelike Jordan Osserman} if the Jordan normal form of $J(x)$
is constant on $S^+(M,g)$ or on $S^-(M,g)$, respectively. Relatively few examples of such
manifolds are known.

The eigenvalue $0$ is distinguished. One says that $(M,g)$ is {\it nilpotent Osserman} if
$J(x)^{p+q}=0$ or equivalently if $0$ is the only eigenvalue of
$J(x)$ for any $x\in TM$. The {\it orders of nilpotency} $n(x)$ and $n(M)$ are then defined by
the properties:
$$J(x)^{n(x)}=0,\quad J(x)^{n(x)-1}\ne0,\quad\text{and}\quad n(M):=\sup_{x\in TM}n(x)\,.$$
Fiedler and Gilkey \cite{refFeGi03} gave examples of $m$ dimensional
pseudo-Riemannian manifolds for any $m\ge4$ where
$n(M)=m-2$; thus
$n(M)$ can be arbitrarily large. However for these examples,
$n(x)$ was constant neither on $S^+(M,g)$ or on $S^-(M,g)$ so these manifolds were neither
spacelike nor timelike Jordan Osserman.

Results of Gilkey and Ivanova \cite{refGiIv} show that if $(M,g)$ is spacelike Jordan Osserman of signature $(p,q)$
where $p<q$, then the Jacobi operator is diagonalizable and hence $(M,g)$ can not be not nilpotent. Thus we suppose
$p\ge q$ henceforth. Examples of spacelike and timelike Jordan Osserman
manifolds of neutral signature $(s,s)$ which are nilpotent of order $2$ have been constructed
Gilkey, Ivanova, and Zhang \cite{refGiIvZh02} for any $s\ge2$. Examples of spacelike
and timelike Jordan Osserman manifolds of signature $(2,2)$ which are nilpotent of order $3$
have been constructed by Garc\'ia-Ri\'o, V\' azquez-Abal and V\' azquez-Lorenzo
\cite{refGVV}. This brief note is devoted to the proof of the following result:
\begin{theorem}\label{thm-1.1}
If $s\ge2$, then there exist pseudo-Riemannian manifolds of signature
$(2s,s)$ which are spacelike Jordan Osserman nilpotent of order $3$ and
which are not timelike Jordan Osserman.
\end{theorem}

Our examples is quite different in flavor from those described in
\cite{refGVV,refGiIvZh02} in several respects. The primary feature is that we are {\bf not} in
the {\it balanced setting} where
$p=q$; the extra timelike directions play a central role in our construction.
Additionally, the examples of \cite{refGVV,refGiIvZh02} are also timelike Jordan Osserman;
this is not the case for our examples.

To prove Theorem \ref{thm-1.1}, it is convenient to work first in a purely algebraic
context. In Section
\ref{Sect2}, we shall construct a family of algebraic curvature tensors $R$ on a vector space
$V$ of signature
$(2s,s)$ which are spacelike Jordan Osserman nilpotent of order $3$ and which are not
timelike Jordan Osserman.  We complete the discussion in Section
\ref{Sect3} by realizing this family geometrically. Our construction will show
that in fact there are many such examples; although we shall use quadratic polynomials to
define the metric in question, this is an inessential feature.

\section{Algebraic curvature tensors}\label{Sect2}

Let $V$ be a finite dimensional real vector space which is equipped with a non-degenerate
symmetric bilinear form
$g(\cdot,\cdot)$ of signature $(p,q)$. Let $R\in\otimes^4V^*$. We say that $R$ is an {\it
algebraic curvature tensor} if $R$ satisfies the symmetries of the Riemann curvature tensor:
\begin{eqnarray}
&&R(x,y,z,w)=-R(y,x,z,w),\nonumber\\
&&R(x,y,z,w)=R(z,w,x,y),\label{eqn-2.a}\\
&&R(x,y,z,w)+R(y,z,x,w)+R(z,x,y,w)=0\,.\nonumber
\end{eqnarray}
The associated Jacobi operator is then defined using equation (\ref{eqn-1.a}) and the notions
spacelike Jordan Osserman and so forth are defined analogously.

Let $s\ge2$. Let $\mathcal{U}:=\{U_1,...,U_s\}$,
$\mathcal{V}:=\{V_1,....,V_s\}$, and $\mathcal{T}:=\{T_1,...,T_s\}$ comprise a basis for
$\mathbb{R}^{3s}$. We let indices $a,b,c,d$ range from $1$ through $s$.

\begin{lemma}\label{lem-2.1} 
Let $g_{ab}=g_{ba}$ be an arbitrary symmetric matrix. Define a metric
$g$ of signature $(2s,s)$ on $\mathbb{R}^{3s}$ whose non-zero components are:
$$g(U_a,U_b)=g_{ab},\quad
g(U_a,V_b)=g(V_b,U_a)=\delta_{ab},\quad g(T_a,T_b)=-\delta_{ab}\,.$$ 
Let $R^{(1)}$ and $R^{(2)}$ be algebraic curvature tensors on
$\Pspan\{U_a\}$. Define a $4$ tensor
$R=R(R^{(1)},R^{(2)})$ on
$\mathbb{R}^{3s}$ whose non-zero entries are
\begin{eqnarray*}
R(U_a,U_b,U_c,U_d)&:=&R^{(1)}(U_a,U_b,U_c,U_d),\\
R(U_a,U_b,U_c,T_d)&=&R(U_a,U_b,T_c,U_d)=R(U_a,T_b,U_c,U_d)\\
  &=&R(T_a,U_b,U_c,U_d):=R^{(2)}(U_a,U_b,U_c,U_d)\,.
\end{eqnarray*}
\begin{enumerate}
\item $R$ is an algebraic curvature tensor on $\mathbb{R}^3$.
\item If $R^{(2)}(U_a,U_b,U_c,U_d):=\delta_{ad}\delta_{bc}-\delta_{ac}\delta_{bd}$, then
$R$ is spacelike Jordan Osserman nilpotent of order $3$ and not timelike Jordan Osserman.
\end{enumerate}
\end{lemma}

\begin{proof} The sum of algebraic curvature tensors is again an algebraic curvature tensor.
If $R^{(2)}=0$, then clearly $R$ is an algebraic curvature tensor since we may assume
$x,y,z,w\in\mathcal{U}$ in establishing the relations of display
(\ref{eqn-2.a}). We may therefore set
$R^{(1)}=0$ and consider only the effect of
$R^{(2)}$ in proving assertion (1). In that case, exactly one of the vectors $x,y,z,w$ must be
taken from $\mathcal{T}$ and the remaining $3$ vectors must be taken from
$\mathcal{U}$. Suppose, for example, $x\in\mathcal{T}$ while
$y,z,w\in\mathcal{U}$. Then replacing $x$ by the corresponding element $\bar
x\in\mathcal{U}$ replaces $R$ by $R^{(2)}$ and thus the relations of display (\ref{eqn-2.a})
follow for $R$ because from the corresponding relations for $R^{(2)}$. This proves assertion
(1).

The tensor $R^{(2)}$ of assertion (2) is the algebraic curvature  tensor of constant sectional
curvature
$+1$ with respect to the standard metric $(U_a,U_b)=\delta_{ab}$. Consequently, it is
invariant under the action of the orthogonal group $O(s)$.

Expand a spacelike vector $X\in\mathbb{R}^{3s}$ in the form $X=u_aU_a+v_aV_a+t_aT_a$ 
where we adopt the Einstein convention and sum over repeated indices.
Then 
$$g(X,X)=g_{ab}u_au_b+2\delta_{ab}u_av_b-\delta_{ab}t_at_b.$$
If $\vec u=0$, then $g(X,X)\le0$. Consequently $\vec u\ne0$. By making an orthogonal
rotation in the $U$ vectors and the same orthogonal rotation in the $V$ and in the $T$ vectors
and by rescaling
$X$, we may suppose without loss of generality that the bases $\mathcal{U}$, $\mathcal{V}$,
and $\mathcal{T}$ have been chosen so that the general form of $g$ and $R$ is the same, so
that $u_1=1$, and so that $u_a=0$ for $a>1$. For $1\le a,b,c,d\le s$, define
$R_{abcd}^{(2)}:=R^{(2)}(U_a,U_b,U_c,U_d)$. Then:
$$\begin{array}{lll}
(J(X)U_a,U_b)=C_{ab},&(J(X)U_a,T_b)=R^{(2)}_{a11b},&(J(X)U_a,V_b)=0,\\
(J(X)T_a,U_b)=R^{(2)}_{a11b},&(J(X)T_a,T_b)=0,&(J(X)T_a,V_b)=0,\\
(J(X)V_a,U_b)=0,&(J(X)V_a,T_b)=0,&(J(X)V_a,V_b)=0,
\end{array}$$
where $C_{ab}=C_{ba}$ is an appropriately chosen matrix. We then have:
$$J(X)U_a=C_{ab}V_b-R^{(2)}_{a11b}T_b,\ 
  J(X)T_a=R^{(2)}_{a11b}V_b,\ J(X)V_a=0\,.$$
It is now clear that $J(X)^3=0$. We have $J(X)X=0$ and $J(X)V_a=0$. Since $R_{a11b}^{(2)}=0$
if
$a=1$ or
$b=1$,
$J(X)T_1=0$. Set $R_{a11b}^{(2)}=\delta_{ab}$ for $a\ge2$. Since $u_1=1$,
$\{X,U_2,...,U_s,T_1,...,T_s,V_1,...,V_s\}$ is a basis for $V$. Consequently:
\begin{eqnarray*}
&&\Range(J(X))=\Pspan\{J(X)X,J(X)U_2,...,J(X)U_s,J(X)T_1,...,J(X)T_s,\\
&&\phantom{\Range(J(X))=\Pspan\{}J(X)V_1,...,J(X)V_s\}\\
&&\phantom{\Range(J(X))}=\Pspan\{J(X)U_2,...,J(X)U_s,J(X)T_2,...,J(X)T_s\}\\
&&\phantom{\Range(J(X))}
  =\Pspan\{C_{2b}V_b-T_2,...,C_{sb}V_b-T_s,V_2,...,V_s\}\,.
\end{eqnarray*}
The set $\{C_{2b}V_b-T_2,...,C_{sb}V_b-T_s,V_2,...,V_s\}$ is linearly independent.
Furthermore:
\begin{eqnarray*}
&&\Range(J(X))\cap\ker(J(X))=\Pspan\{V_2,...,V_s\},\\
&&\Range(J(X)^2)=\Pspan\{V_2,...,V_s\}\,.
\end{eqnarray*}
It is now clear that $R$ is spacelike Jordan Osserman nilpotent of order $3$.
Since $J(T_1)=0$ while $J(U_1-V_1)=J(U_1)\ne0$, $R$ is not timelike Jordan
Osserman.
\end{proof}

\section{Geometric Realizations}\label{Sect3}

We complete the proof of Theorem \ref{thm-1.1} by showing that the structures of Lemma
\ref{lem-2.1} are geometrically realizable. The metrics we shall consider are similar those
described in different contexts in
\cite{Br00,Ka00,refOpr}. We take coordinates of the form
$(u_1,...,u_s,v_1,...,v_s,t_1,...,t_s)$ on
$\mathbb{R}^{3s}$. Let
$$U_a:=\textstyle\frac{\partial}{\partial u_a},\quad
  V_a:=\textstyle\frac{\partial}{\partial v_a},\quad\text{and}\quad
  T_a:=\textstyle\frac{\partial}{\partial t_a}$$
be the associated coordinate frame for the tangent bundle. We let the index $r$ range from
$1$ to $3s$ and index the full coordinate frame 
$$\{e_1,...,e_{3s}\}:=\{U_1,...,U_s,V_1,...,V_s,T_1,...,T_s\}\,.$$
Theorem \ref{thm-1.1} will follow from Lemma \ref{lem-2.1} and from the following Lemma:

\begin{lemma}\label{lem-3.1} Let $R^{(2)}$ be a fixed algebraic curvature tensor on
$\mathbb{R}^s$.
 Define a metric $g$ of signature
$(2s,s)$ on $\mathbb{R}^{2s,s}$ whose non-zero inner products are given by:
\begin{eqnarray*}
&&g(U_a,U_b)=\psi_{abcd}u_ct_d\text{ where }\psi_{abcd}=\psi_{bacd}:=\textstyle
   -\frac23(R_{acdb}^{(2)}+R_{adcb}^{(2)}),\\
&&g(U_a,V_b)=g(V_b,U_a)=\delta_{ab},\quad\text{and}\quad
g(T_a,T_b)=-\delta_{ab}\,.
\end{eqnarray*}
Let $R^{(1)}_{abcd}(u,t):=R(U_a,U_b,U_c,U_d)(u,t)$. Then
$R(u,t)=R(R^{(1)}(u,t),R^{(2)})$.
\end{lemma}

\begin{proof} At this point, we change our indexing convention slightly for the remainder of
the proof. We shall let indices
$a,b,c$ index elements of
$\mathcal{U}$, indices $\alpha,\beta,\gamma$ index elements of $\mathcal{V}$, and indices
$i,j,k$ index elements of $\mathcal{T}$. Indices $r_\mu$ will index the full coordinate
basis. By an abuse of notation, we shall set
$\Gamma_{abc}=g(\nabla_{U_a}U_b,U_c)$,
$\Gamma_{abi}=g(\nabla_{U_a}U_b,T_i)$, etc.
We replace an element of $\mathcal{T}$ by the corresponding element of $\mathcal{U}$ to define
$\tilde\psi_{abci}$,
$\tilde R^{(2)}_{abci}$,
$\tilde R^{(2)}_{abic}$, $\tilde R^{(2)}_{aibc}$, and $\tilde R^{(2)}_{iabc}$. The non-zero
Christoffel symbols of the metric are:
\begin{equation}\label{eqn-3.a}
\begin{array}{l}
\Gamma_{abc}=\textstyle\frac12(\tilde\psi_{bcai}+\tilde\psi_{acbi}-\tilde\psi_{abci})t_i,\\ 
\Gamma_{iab}=\Gamma_{aib}=-\Gamma_{abi}=\textstyle\frac12\tilde\psi_{abci}u_c\,.
 \vphantom{\vrule height 14pt}
\end{array}\end{equation}
We raise indices to see:
\begin{equation}\label{eqn-3.b}
\Gamma_{r_1r_2}{}^a=0,\quad\Gamma_{r_1r_2}{}^i=-\Gamma_{r_1r_2i},
\quad\text{and}\quad\Gamma_{r_1r_2}{}^\alpha=\Gamma_{r_1r_2a}\,.\end{equation}
The curvature tensor is given by:
$$
R_{r_1r_2r_3r_4}=e_{r_1}\Gamma_{r_2r_3r_4}-e_{r_2}\Gamma_{r_1r_3r_4}
  +\Gamma_{r_1r_5r_4}\Gamma_{r_2r_3}{}^{r_5}-\Gamma_{r_2r_5r_4}\Gamma_{r_1r_3}{}^{r_5}\,.
$$
If $r_5$ indexes an element of $\mathcal{V}$, then $\Gamma_{\star r_5\star}=0$ by equation
(\ref{eqn-3.a}) while if $r_5$ indexes an element of $\mathcal{U}$, then
$\Gamma_{\star\star}{}^{r_5}=0$ by equation (\ref{eqn-3.b}). Thus $r_5$ must index an element
of $\mathcal{T}$ and consequently, we may express:
\begin{equation}\label{eqn-3.c}
R_{r_1r_2r_3r_4}=e_{r_1}\Gamma_{r_2r_3r_4}-e_{r_2}\Gamma_{r_1r_3r_4}
  +\Gamma_{r_1ir_4}\Gamma_{r_2r_3}{}^{i}-\Gamma_{r_2ir_4}\Gamma_{r_1r_3}{}^{i}\,.
\end{equation}
Thus by equation (\ref{eqn-3.a}), quadratic terms in $\Gamma$ can only appear in equation
(\ref{eqn-3.c}) if
$r_1$, $r_2$, $r_3$, and $r_4$ all index elements of $\mathcal{U}$.
The only other non-zero curvatures occur when exactly one of $r_\nu$ indexes an
element of $\mathcal{T}$ and the remaining $r_\nu$ index elements of $\mathcal{U}$. We may
therefore compute the proof by computing:
\medbreak\qquad$R(U_a,U_b,U_c,T_i)=U_a\Gamma_{bci}-U_b\Gamma_{aci}=
\textstyle\frac12(\tilde\psi_{acbi}-\tilde\psi_{bcai})$
\smallbreak\qquad\qquad\qquad
$=-\textstyle\frac13(\tilde R_{abic}^{(2)}+
   \tilde R_{aibc}^{(2)}-\tilde R_{baic}^{(2)}-\tilde R_{biac}^{(2)})$
\smallbreak\qquad\qquad\qquad
$=-\textstyle\frac13(2\tilde R^{(2)}_{abic}-\tilde R^{(2)}_{iabc}-\tilde R^{(2)}_{ibca})$
$=-\textstyle\frac13(2\tilde R^{(2)}_{abic}+\tilde R^{(2)}_{icab})=\tilde R^{(2)}_{abci}$.
\end{proof}

\begin{remark}\label{rem-3.2}
\rm It is worth giving a very specific example. Define an inner product $g$ on
$\mathbb{R}^6$ whose non-zero components are, up to the usual $\mathbb{Z}_2$ symmetries given
by:
\begin{eqnarray*}
&&g(U_1,U_1)=-2u_2t_2,\quad g(U_2,U_2)=-2u_1t_1,\quad g(U_1,U_2)=u_1u_2,\\
&&g(U_1,V_1)=g(U_2,V_2)=-g(T_1,T_1)=-g(T_2,T_2)=1.\,\end{eqnarray*}
This manifold has signature $(4,2)$. It is spacelike Jordan Osserman nilpotent of order $3$.
It is not timelike Jordan Osserman. Furthermore, it is curvature homogeneous up to order $0$ as defined by
Kowalski, Tricerri, and Vanhecke
\cite{KTV92}.
\end{remark}

\section*{Acknowledgments} Research of P. Gilkey partially supported by the NSF (USA) and
MPI (Leipzig). Research of S. Nik\v cevi\'c partially supported by the Dierks Von Zweck
Stiftung (Essen), DAAD (Germany) and MM 1646 (Srbija). It is a pleasant task to thank Professor E.
Garc\'{\i}a--R\'{\i}o for helpful discussions on this matter.

\end{document}